\documentclass[12pt]{article}%
\usepackage[applemac]{inputenc}
\usepackage{amsmath,amssymb}
\usepackage{amsfonts}
\usepackage{amsmath,amssymb}
\usepackage[applemac]{inputenc}
\usepackage{amsmath,amssymb}
\usepackage{color}
\usepackage{color, colortbl}
\usepackage{amsmath}
\usepackage{amssymb}
\usepackage{graphicx}
\usepackage{tikz}
\usepackage{geometry}
\usetikzlibrary{arrows}
\usepackage{amsfonts}
\usepackage{amsmath,amssymb}
\usepackage[applemac]{inputenc}
\usepackage{amsmath,amssymb}
\usepackage{color}
\usepackage{amsmath}
\usepackage{amssymb}
\usepackage{graphicx}
\usepackage{tikz}
\newtheorem{theorem}{Theorem}[section]
\newtheorem{lemma}[theorem]{Lemma}
\newcommand{\E}{E}

\newtheorem{remark}[theorem]{Remark}

\usetikzlibrary{arrows,shapes,automata,petri}
\newenvironment{myenumerate}{

\begin{enumerate}}{\end{enumerate}}
\newcommand{\dproof}{\noindent {Proof.} \quad}
\newcommand{\fproof}{\hfill $\square$ \bigskip}

\numberwithin{equation}{section}

\definecolor{LightCyan}{rgb}{0.88,1,1}

\def\E{\mathbb{E}}

\def\1B{\text{1\!\!I}}

\def\<{\langle}
\def\>{\rangle}

\def\R{\mathbb{R}}

\def\U{\mathcal{U}}
\def\P{\mathbb{P}}
\def\Y{\widehat{Y}}

\def\L{\mathcal{L}}
\def\z{\zeta}

\begin{document}

\title{A Kalman filter for linear systems driven by time-space Brownian sheet}
\author{Nacira Agram$^{1},$ Bernt \O ksendal$^{2}$, Frank Proske$^{2}$ and Olena Tymoshenko$^{2,3}$}
\date{8 July 2024\\
\vskip 0.7cm
In memory of Habib Ouerdiane}
\maketitle

\footnotetext[1]{Department of Mathematics, KTH Royal Institute of Technology 100 44, Stockholm, Sweden. \newline
Email: nacira@kth.se. }

\footnotetext[2]{%
Department of Mathematics, University of Oslo, Norway. \\
Emails: oksendal@math.uio.no, proske@math.uio.no, otymoshenkokpi@gmail.com}
\footnotetext[3]{%
Department of Mathematical Analysis and Probability Theory,Igor Sikorsky Kyiv Polytechnic Institute, Ukraine. \\
}

\begin{abstract}
We study a linear filtering problem where the signal and observation processes are described as solutions of linear stochastic differential equations driven by time-space Brownian sheets. We derive a stochastic integral equation for the conditional value of the signal given the observation, which can be considered a time-space analogue of the classical Kalman filter. The result is illustrated with examples of the filtering problem involving noisy observations.
\end{abstract}

\textbf{Keywords :} Filtering, linear stochastic differential equations, time-space Brownian sheets, conditional expectation, Riccati equation, Kalman type filter.

\textbf{MSC 2020 :} 60G15, 60G35, 60G60, 60H15, 60H20, 62M20, 93E10, 93E11, 94AXX

\section{Introduction}
The Kalman filter, introduced in the 1960s by R.E. Kalman \cite{K}, revolutionized the field of signal processing and control theory by providing an efficient recursive solution to the linear quadratic estimation problem. The Kalman filter consists of a series of mathematical equations that offer an effective  recursive method for estimating the state of a process while minimizing the mean squared error. Traditional applications of the Kalman filter typically involve time-evolving signals and observations influenced by temporal Brownian motion. For more details about linear filtering, we refer to \cite{J,Oksendal13, Walsh}, and for nonlinear filtering, we refer to \cite{CR}.

However, many modern applications in fields such as spatial statistics, environmental science, and image processing require dealing with systems influenced by both spatial and temporal variability. As the Kalman filter has quickly become an important component of modern control systems theory and practice, it is natural to expand the concepts of the Kalman filtering into scientific fields where the methodology can be adapted to solve numerous state-space oriented problems. Moreover, Kalman filtering (of spatial and temporal data) is still an active area of research in atmospheric and oceanic sciences (see \cite{Can}, \cite{Mil}). 

In this paper, we extend the classical Kalman filter framework to encompass systems where both the signal and observation processes are driven by time-space Brownian sheets. The time-space Brownian sheet, a natural generalization of the one-dimensional Brownian motion to higher dimensions, provides a more comprehensive model for phenomena that evolve over both time and space. We believe that this extension, incorporating the type of noise present in the system, could enable new applications and yield better estimates than the classical theory in certain scenarios.

Our work focuses on deriving a time-space analogue of the Kalman filter that retains the recursive estimation property of its classical counterpart. 

The organization of the paper is as follows:\\
We begin, in Section 2, by formulating the linear time-space filtering problem within the framework of stochastic partial differential equations (SPDEs), describing the dynamics of the signal and observation processes driven by time-space Brownian sheets. Subsequently, in Section 3, we give a solution formula for the signal process and an $L^2$ estimate of it. Section 4 is devoted to the equation given by the innovation process. In Section 5, we provide a representation formula for the estimate, commonly known as the Kalman filter, within the time-space framework. Both the differential and integral forms are discussed in Section 6.
Finally, we apply our results to the filtering problem involving noisy observations of a given random variable, and noisy observation of a Brownian sheet.

\section{The linear time-space filtering problem}
In this section, we aim to formulate the time-space filtering problem. Therefore, we shall describe the dynamics of the signal process, which evolves according to an SPDE driven by a time-space Brownian sheet. Similarly, we model the observation process, which is also influenced by a time-space Brownian sheet. 

To provide context, let us first recall the situation in the classical (single-parameter) Brownian motion case (e.g. see \cite{Oksendal13}):\\
Suppose the signal process $X(t)$  is described by a 1-dimensional stochastic differential equation of the form
\begin{align*}
dX(t)=F(t) X(t) dt + C(t) dB_1(t); \quad X(0)=X_0,
\end{align*}
and the observation process is given by
\begin{align*}
dU(t)=G(t) X(t) dt + D(t) dB_2(t);  \quad U(0)=\E[X_0].
\end{align*}
Here $B_1(t),B_2(t)$ are independent Brownian motions, $X_0$ is a given Gaussian random variable, independent of $(B_1,B_2)$ and $F(t),C(t),G(t),D(t) \in \R$ are 
given bounded (Borel measurable) deterministic functions, $D(t)$ is bounded away from 0.\\
Let $\U=\{\U_t\}_{t\geq 0}$ denote the filtration generated by the observation process. The problem is to find the best estimate of the signal at time $t$, given the observation up to time $t$. Mathematically this means that we want to find the conditional expectation
$$\widehat{X}(t):=\E[X(t)| \U_t],$$
where $\E$ denotes expectation with respect to the probability law $\P$ of the 2-dimensional Brownian motion $(B_1,B_2)$.
This conditional expectation  $\widehat{X}(t)$, called \emph{the Kalman filter}, is given by the SDE
\begin{align}
d\widehat{X}(t)=\Big( F(t) - \frac{G^2(t)S(t)}{D^2(t)} \Big)\widehat{X}(t)dt + \frac{G(t)S(t)}{D^2(t)} dU(t); \quad \widehat{X}(0)=\E[X_0], \label{kalman1}
\end{align}
where $S(t):=\E[(X(t)-\widehat{X}(t))^2]$ is the error process. It satisfies the (deterministic) Riccati equation
\begin{align}
\frac{d}{dt}S(t)=2 F(t) S(t) - \frac{G^2(t)}{D^2(t)} S^2(t) +C^2(t), \quad
S(0)=\E[(X(0)-\E[X(0)])^2]. \label{riccatti1}
\end{align}

Note that, using \eqref{riccatti1}, $S(t)$ can be computed beforehand, and then the Kalman equation \eqref{kalman1} allows us to update the estimate as the observation progresses. 

\vspace{1cm}

We now proceed to the time-space (Brownian sheet) extension of this system. Recall, that a one-dimensional Brownian sheet is a 2-parameter, centered Gaussian process $B = {B(t,x); t,x \geq 0}$ whose covariance is given by
$$\E[B(t,x)B(t',x')]=\min(t,t')\times \min(x,x'), \,\,\textnormal{for all}\,\,\, t,t',x,x' \geq 0.$$
Note that in the 1-parameter case we have the classical Brownian motion
$B(t)$. 

Let $(B_1,B_2)=(B_1(t,x),B_2(t,x))$ be a 2-dimensional Brownian sheet defined on a filtered probability space  $(\Omega,\mathcal{F},\mathbb{F}=\{\mathcal{F}_t\}_{t\geq 0},\P)$. 
Suppose the signal process $Y(z)$ is given by the SPDE
\begin{align}
dY(z)=F(z) Y(z) dz + C(z) B_1(dz); z=(t,x)\geq 0; \quad Y(0,x))=Y(t,0)=Y_0 \in \R, \label{signal1}
\end{align}
and the observation process $U(z)$ is given by the SPDE
\begin{align}
dU(z)=G(z) Y(z) dz + D(z) B_2(dz); z=(t,x) \geq 0; \quad U(0,x)=U(t,0)=U_0 \in \R. \label{obs1}
\end{align}
Here $z=(t,x)$ and $B_1(t,x),B_2(t,x)$ are independent Brownian sheets, and $F(z)$,$C(z)$, $G(z)$ and $D(z)$ are 
given bounded (Borel measurable) deterministic functions, $D(z)$ is bounded away from 0. We assume that $Y_0$ is a given Gaussian random variable, independent pf $(B_1,B_2)$.

The differential equations \eqref{signal1}, \eqref{obs1} are short hand notations  for the following SPDEs (with the same boundary conditions):
\begin{align}
&\text{(signal)}\quad \frac{\partial ^2}{\partial t \partial x}Y(t,x)= F(t,x)Y(t,x) +C(t,x) \overset{\bullet}{B}_1(t,x)\label{signal2}\\
&\text{(observation)}\quad \frac{\partial ^2}{\partial t \partial x}U(t,x)= G(t,x)Y(t,x) +D(t,x) \overset{\bullet}{B}_2(t,x)\label{signal2.2},
\end{align}
where, in the sense of distributions,  $\overset{\bullet}{B}_j(t,x)=\frac{\partial ^2}{\partial t \partial x}{B}_j(t,x); j=1,2.$ \\

Let $\U=\{\U_z\}_{z\geq 0}$ denote the filtration generated by the observation process $U(z)$, i.e. $\mathcal{U}_z$ is the sigma-algebra generated by the random variables $\{ U(s,a)\}_{s\leq t, a\leq x}$.\\
The problem is to find the best estimate of the signal at time-space $z=(t,x)$, given the observations up to time-space  $(t,x)$. We aim to determine the conditional expectation of $Y(z)$ with respect to the $\sigma$-algebra  $\mathcal{U}_z$, denoted as
$$\Y(z):=\E[Y(z)| \U_z],$$
where $\E$ represents the expectation  based on the  probability law $\P$ of the 2-dimensional Brownian sheet $(B_1,B_2)$.\\

To study this problem we will follow the approach in Ch.6 in \cite{Oksendal13}, but with necessary (and nontrivial) modifications. 
We first consider some auxiliary results:

\subsection{Relation to the projection operator}

For estimation of the signal process $Y(z)$, given by (\ref{signal1}),  we  prove the following statement.
\begin{lemma}
    The processes $Y(z)$ and $U(z)$ are Gaussian.
\end{lemma}
\dproof We first prove that $Y(z)$ is Gaussian:\\
Given $Y(0)=Y_0$ (Gaussian) we define inductively
\begin{align*}
    Y_k(z)= Y_0 + \int_0^z f(\zeta) Y_{k-1}(\zeta) d\zeta + \int_0^z C(\zeta) B_1(d\zeta);\,\, k=1,2, ...
\end{align*}
Then by induction we see that $Y_k(z)$ is Gaussian for all $k$. Moreover, it is easy to see that $Y_k(z) \rightarrow$ $Y(z)$ in $L^2(\P)$ as $k \rightarrow \infty$. Hence $Y(z)$ is Gaussian.\\
The same argument proves that $U(z)$ is Gaussian for all $z$.
\fproof

For fixed $Z=(T,X)$ we let $\L(U)=\L(U,Z)$ denote the closure in $L^2(\P)$ of the set of all linear combinations of the form
\begin{align*}
    c_0 + c_1 U(z_1) + c_2 U(z_2) + ... + c_k U(z_k),
\end{align*}
where $c_j\in \R$ are constants and $z_j=(t_j,x_j) \leq Z$ (i.e. $t_j \leq T,x_j \leq X$) for all $j$.\\
Let
$$\mathcal{P}_{\L}(\cdot): L^2(\P) \mapsto \L(U)$$
denote the orthogonal projection from $L^2(\P)$ down to $\L(U)$.
The following lemma describes the form of the best estimate for the process $Y(z)$ from the equation (\ref{signal2}).
\begin{lemma}
The best estimate of $Y(z)$ coincides with the projection of $Y(z)$ down to $\L(U,z)$, i.e.
$$\Y(z)=\mathcal{P}_{\L(U,z)}(Y(z)).$$
\end{lemma}
\dproof  Define $\overset{\vee}{Y}(z)= \mathcal{P}_{\L(U,z)}(Y(z))= \mathcal{P}_{\L}(Y(z))$. We claim that $Y(z)-\overset{\vee}{Y}(z)$ is independent of $\U_z$.
To see this, note that $\mathcal{P}_{\L(U,z)}(Y(z))$ is Gaussian and hence 
$$ (\overset{\vee}{Y}(z), U(z_1), U(z_2), ... ,U(z_n))$$
is Gaussian for all $z_j$, $j=1,2,...$ Moreover, $Y(z) - \overset{\vee}{Y}(z)$ is orthogonal to $\L(U,z)$. i.e.
$$\E[(Y(z) - \overset{\vee}{Y}(z)) U(z_j)]=0 \text { for all } z_j.$$
Hence $Y(z)-\overset{\vee}{Y}(z)$ and $U(z_j)$ are uncorrelated. Since they are Gaussian they are also independent. Hence $Y(z)-\overset{\vee}{Y}(z)$ is independent of $\U_z$. Therefore
$$\E[\chi_{H} (Y(z)-\overset{\vee}{Y}(z))]=\P(H) \E[Y(z)-\overset{\vee}{Y}(z)]=0\,\,\, \textnormal{ for all}\,\,\, H \in \U_z.$$ It follows that
$$0=\E[\chi_{H}(Y(z)-\mathcal{P}_{\L}(Y(z)))]=\E[\chi_{H} Y(z)]- \E[\chi_{H} \mathcal{P}_{\L}(Y)].$$
Since this holds for all $H \in \U_z$ we can conclude that $\mathcal{P}_{\L}(Y(z))= \E[Y(z) | \U_z]= \widehat{Y}(z).
$
\fproof

\section{The signal process}
We begin this section by outlining the formulation of the signal process as governed by an SPDE driven by a time-space Brownian sheet. Following this, we provide a rigorous derivation of the solution formula. Additionally, we establish an $L^2$ estimate to quantify the mean-square stability and accuracy of the solution.

We first prove the following result:
\begin{lemma}\label{pde0}
\begin{myenumerate}
\item
Let $h(t,x)$ be a given bounded function.    The solution $y(t,x)$ of the PDE
    \begin{align}\label{pde}
    \frac{\partial^2}{\partial t \partial x}y(t,x)=h(t,x) y(t,x); (t,x) \in \R_{+}\times \R_{+}; \quad y(0,x)=y(t,0)=y_0
    \end{align}
    is given by
\begin{align}\label{pde2}
&y(t,x)=y_0\sum_{n=0}^{\infty} J_n h(t,x),
\end{align}
where 
\begin{align}
J_n h(t,x):= \int_0^{(t,x)} \big( \int_0^{u_1} ... \int_0^{u_n} h(u_1) h(u_2) ... h(u_{n+1}) du_1 du_2 ...du_n \big) du_{n+1}.
\end{align}
\item 
For given $\zeta \geq 0$ the solution $Y(t,x)$ of the signal equation
\begin{align*}
    dY(z)=F(z) Y(z)dz + C(z) B_1(dz); z\geq \zeta; \quad Y(\zeta) \text{ given}
\end{align*}
has the form
\begin{align}\label{2.10}
    Y(z)=& Y(\zeta)\sum_{n=0}^{\infty}J_n F(\zeta;z) \nonumber\\
   &+ \text{ terms involving } dB_1\text{ -integrals from } \zeta \text { to } z
   \end{align}
   where
   \begin{align}
   &J_n F(\zeta;z):= \int_{\zeta}^{z} \big( \int_{\zeta}^{u_1} ... \int_{\zeta}^{u_n} F(u_1) F(u_2) ... F(u_{n+1}) du_1 du_2 ...du_n \big) du_{n+1}.\label{Jn}
\end{align}
\end{myenumerate}   
\end{lemma}
\dproof
(i) The integral version of \eqref{pde} is, putting $z=(t,x)$, 
\begin{align*}
    y(z)=y_0+\int_0^z h(\zeta)y(\zeta)d\zeta; \quad y(0,x)=y(t,0)=y_0.
\end{align*}
Substituting $y(\zeta)=y_0+\int_0^{\z} h(\zeta_1) y(\zeta_1)d\zeta_1$ in the above and repeating by induction we get
\begin{align*}
y(z)= y_0 \sum_{n=0}^\infty \int_0^{z} \big( \int_0^{u_1} ... \int_0^{u_n} h(u_1) h(u_2) ... h(u_{n+1}) du_1 du_2 ...du_n \big) du_{n+1}.\
\end{align*}\\
(ii) This follows by the same argument as in part (i).
\fproof

\begin{remark}
    Note that if the function $F$ has the form $F(t,x)=f(t)g(x)$, then \eqref{Jn} simplifies to
    \begin{align*}
J_n F((s,a);z)=\frac{1}{n! n!}(\int_{s}^t f(u)du)^n (\int_{a}^x g(v)dv)^n;\quad z=(t,x). 
\end{align*}
\end{remark}
We now apply this to the signal process $Y(z)$ given by the integral equation
\begin{align*}
    Y(z)=Y(0)+\int_0^z F(\zeta)Y(\zeta)d\zeta + \int_0^z C(\zeta) B_1(d\z);\quad Y(0)=Y_0 \in \R.
\end{align*}
Taking expectation we get
\begin{align*}
    \E[Y(z)]=\E[Y(0)]+\int_0^z F(\zeta) \E [Y(\zeta)]d\zeta.
\end{align*}
Applying Lemma \ref{pde0} to $y(z):=\E[Y(z)]$ we get
\begin{align}
\E [Y(z)]= \E[Y(0)] \sum_{n=0}^\infty \int_{0}^{z} \big( \int_{0}^{u_1} ... \int_{0}^{u_n} F(u_1) F(u_2) ... F(u_{n+1}) du_1 du_2 ...du_n \big) du_{n+1}.\label{2.12a}
\end{align}
We also need an expression for $\E[Y^2(z)]$. To find this expression, we use the It\^o formula for systems driven by the Brownian sheet \cite {WZ}. So, we have
\begin{align*}
\E[Y^2(z)]&=\E[Y_0]^2 + \E\Big[\int_0^z \Big\{ 2Y(\z)F(\z) Y(\z) + C^2(\z)\Big\} d\z\nonumber\\
&+ \int_0^z \int_0^z I(\z \bar{\wedge}\z')\Big\{F(\z)Y(\z)Y(\z') + F(\z')Y(\z')Y(\z)\Big\} d\z d\z' \Big]\nonumber\\
&=Y_0^2 + \int_0^z \Big\{ 2F(\z)\E\Big[Y^2(\z) \Big]+C^2(\z)\Big\} d\z\nonumber\\
&+ \int_0^z \int_0^z I(\z \bar{\wedge}\z')\Big\{F(\z)\E\Big[Y(\z)Y(\z')\Big] 
+ F(\z')\E\Big[Y(\z')Y(\z) \Big]\Big\} d\z d\z'.\nonumber\\
\end{align*}

By Lemma 5.1 in \cite{AOPT2}, we have
\begin{align}
\frac{\partial ^2}{\partial t \partial x} \E\Big[Y^2(z)\Big]
&=2F(z)\E\Big[Y^2(z)\Big]+ C^2(z)\nonumber\\
&+\E\Big[\Big(\int_0^tF(\z_1,x) Y(\z_1,x)d\z_1\Big) \Big(\int_0^x F(t,\z_2')Y(t,\z_2') d\z_2'\Big)\Big]\nonumber\\
&=2F(z)\E\Big[Y^2(z)\Big] +C^2(z)\nonumber\\
&+\int_0^t \int_0^x F(\z_1,x)F(t,\z_2') \E\Big[Y(\z_1,x)Y(t,\z_2')\Big] d\z_1 d\z_2' \nonumber\\
&=2F(z)\E\Big[Y^2(z)\Big] +C^2(z)\nonumber\\
&+\int_0^t \int_0^x F(s,x)F(t,a) \E\Big[Y^2(s,a)\Big] ds da.\label{2.14a}
\end{align}

\section{The innovation process}
We first state a useful observation about the linear span of the process $U(\z); \z \in [0,T] \times [0,X]$:
\begin{lemma}
    Put $Z=(T,X)$. Then
    \begin{align}
\L(U,Z)= \Big\{ c_0 + \int_0^Z f(\z)U(d\z);\quad c_0 \in \R, f \in L^2([0,T]\times [0,X]) \text{ (deterministic)} \Big\}.
    \end{align}
\end{lemma}

\dproof
This follows from the definition of $\L(U,Z)$.
\fproof

Now we define the innovation process $N(z)$ as follows:
\begin{equation}\label{4.21}
  N(z)= U(z)-\int_0^z G(\z)\widehat{Y}(\z) d\z;  \textnormal{ i.e.}
\end{equation}
 \begin{equation} \label{4.3}
 N(dz)=G(z)\Big(Y(z)-\widehat{Y}(z)\Big) dz +D(z) U(dz).
\end{equation}

\begin{remark}
    We mention that one can find an $\mathcal{U}$-predictable version of $\widehat{Y}$ by using predictable projections of measurable 2-parameter processes. (See P. Imkeller \cite{I2}.)
\end{remark}
Next lemma present properties of the process $N(z)$. 
\begin{lemma} Let process $N(z)$ be defined as in (\ref{4.21})((\ref{4.3})). Then
\begin{myenumerate}
    \item
    $N(z)$ is Gaussian for all $z$;
    \item
    $\E[N(z)]=0$ for all $z$;
    \item
    $\E[N^2(z)]= \int_0^z D^2(\z) d\z $;
    \item
    $N$ has orthogonal (and hence independent) increments.
    \item
    Define
    \begin{align}
        M(dz)&=\frac{1}{D(z)} N(dz)= \frac{1}{D(z)}U(dz) -\frac{G(z)}{D(z)} \widehat{Y}(z)dz \nonumber\\
        &=\frac{G(z)}{D(z)} \big[Y(z)-\widehat{Y}(z)\big]dz + B_2(dz).\label{M}
    \end{align}
 Then
 $\E[M(s,a) M(s',a')]= \min(s,s') \min(a,a')$; \quad  for all $(s,a),(s',a')$.
 \item 
 $M(z)$ is a Brownian sheet.   
    \end{myenumerate} 
\end{lemma}
\dproof These results follow easily from the definition of $N$ and $M$. We skip the details.
\fproof
 \section{A representation formula for $\widehat{Y}(z)$}
In this section we prove the following:
 \begin{lemma}\label{4.1}
 For all $z$ the following holds:
 \begin{align*}
     \widehat{Y}(z) = \E[Y(z)] + \int_0^z \frac{\partial^2}{\partial s \partial a} \E[Y(z) M(s,a)]M(ds,da).
     \end{align*}
 \end{lemma}
 \dproof Since $\L(M)=\L(U)$ we can for each $z$ find a function $g(z,\z)$ such that
 \begin{align*}
     \widehat{Y}(z)= c_0(z)+ \int_0^z g(z,\z) M(d\z),
 \end{align*}
 where $c_0(z)=\E[\widehat{Y}(z)]=\E[Y(z)].$
Moreover,
\begin{align*}
    Y(z)-\widehat{Y}(z) \perp \int_0^z f(\z)M(d\z)
\end{align*}
for all $f \in L^2([0,T]\times [0,X]).$
Therefore, by the Ito isometry,
\begin{align*}
\E[Y(z) \int_0^z &f(\z) M(d\z)]= \E[\widehat{Y}(z) \int_0^z f(\z) M(d\z)]\nonumber\\
&=\E\Big[\Big(\int_0^z g(z,\z) M(d\z)\Big) \Big( \int_0^z f(\z) M(d\z)\Big)\Big]=\int_0^z g(z,\z) f(\z) d\z.
    \end{align*}
In particular, choosing
$f(\z)=\chi_{[0,s]\times [0,a]}(\z)$, 
where $z=(t,x), s\leq t, a\leq x,$ we get
\begin{align*}
    \E[Y(z)M(s,a)]= \int_0^s \int_0^a g(z,\z) d\z,
\end{align*}
and hence
\begin{align*}
    g(z,(s,a))=\frac{\partial^2}{\partial s \partial a} \E[Y(z) M(s,a)].
\end{align*}
\fproof

\section{A stochastic integral equation for $\widehat{Y}(z)$}
We will discuss both the differential and integral forms of the Kalman filter.\\
From \eqref{M} we have, with $\widetilde{Y}=Y - \widehat{Y},$
\begin{align*}
    M(s,a)=\int_0^{(s,a)}\frac{G(\z)}{D(\z)} \widetilde{Y}(\z)d\z + B_2(s,a).
\end{align*}
Using this, and that $Y(z)$ is independent of $B_2(\cdot)$, we get, for $s\leq t,a\leq x, z=(t,z),$
\begin{align*}
    \E[Y(z)M(s,a)]=\int_0^{(s,a)}\frac{G(\z)}{D(\z)} \E[Y(z)\widetilde{Y}(\z)] d\z.
\end{align*}
For $\z \leq z$ we have, by \eqref{2.10},
\begin{align*}
    &\E\big[Y(z)\widetilde{Y}(\z)\Big]\\
    &= \E\Big[\Big(Y(\zeta) \sum_{n=0}^{\infty}J_n F(\zeta;z)+ \text{ terms involving } dB_1\text{ -integrals from $\zeta$ to $z$}\Big)\widetilde{Y}(\z)\Big]\\
    &=\sum_{n=0}^{\infty}J_n F(\zeta;z)E\Big[Y(\zeta) \widetilde{Y}(\z)\Big]\\
    &=\sum_{n=0}^{\infty}J_n F(\zeta;z)\E\Big[(Y(\zeta)-\widehat{Y}(\zeta)) \widetilde{Y}(\z)\Big]\\
    &=\sum_{n=0}^{\infty}J_n F(\zeta;z)\E\Big[\widetilde{Y}(\z)^2\Big]
\end{align*}
This gives
\begin{align*}
    \E[Y(z)M(s,a)]&= \E\Big[Y(z)\Big(\int_0^{(s,a)}\frac{G(\z)}{D(\z)}(Y(\z)-\widehat{Y}(\z))d\z +B_2(s,a)\Big)\Big]\nonumber\\
    &=\int_0^{(s,a)}\frac{G(\z)}{D(\z)}\E[Y(z)\widetilde{Y}(\z)] d\z\nonumber\\
    &=\int_0^{(s,a)}\frac{G(\z)}{D(\z)}\sum_{n=0}^{\infty}J_n F(\zeta;z)\E\Big[\widetilde{Y}(\z)^2\Big]d\zeta.
\end{align*}
Hence
\begin{align}
  \frac{\partial^2}{\partial s \partial a} \E[Y(z) M(s,a)]=\frac{G(s,a)}{D(s,a)}\sum_{n=0}^{\infty}J_n F((s,a);z)\E\Big[\widetilde{Y}(s,a)^2\Big].\label{5.51} 
\end{align}
From Lemma \ref{4.1} we have
\begin{align}
     \widehat{Y}(z) = \E[Y(z)] + \int_0^z \frac{\partial^2}{\partial s \partial a} \E[Y(z) M(s,a)]M(ds,da).\label{5.6}
     \end{align}
 Combining this with \eqref{5.51} we get
 \begin{align}
     \widehat{Y}(z) = \E[Y(z)] + \int_0^z \Big(\frac{G(s,a)}{D(s,a)} \sum_{n=0}^{\infty}J_n F((s,a);z)\E\Big[\widetilde{Y}(s,a)^2\Big]\label{5.3}\Big)  M(ds,da).     
     \end{align}
Recalling that
     \begin{align*}
         M(dz)&=\frac{1}{D(z)} N(dz)= \frac{1}{D(z)}U(dz) -\frac{G(z)}{D(z)} \widehat{Y}(z)dz
     \end{align*}
we obtain the first part of the following theorem, which is our main result:\\
\begin{theorem}{(The time-space Kalman filter)} \label{5.1}
\begin{myenumerate}
    \item 
The best estimate 
$$\widehat{Y}(z)=\E[Y(z)|\mathcal{U}_z]$$ 
of $Y(z)$ given the observations $\{U(\z);\z \leq z \}$ satisfies the following stochastic integral equation:
\small
     \begin{align}
     &\widehat{Y}(z) 
     = \E[Y(z)] - \int_0^z \frac{G^2(s,a)}{D^2(s,a)}\sum_{n=0}^{\infty}J_n F((s,a);z)\E\Big[\widetilde{Y}(s,a)^2\Big]\widehat{Y}(s,a)dsda\nonumber\\
     &+\int_0^z \frac{G(s,a)}{D^2(s,a)}\sum_{n=0}^{\infty}J_n F((s,a);z) \E\Big[\widetilde{Y}(s,a)^2\Big]U(ds da)\label{5.2}      
     \end{align}  
     \item
In differential form this can be written
     \begin{align*}
 d\widehat{Y}(z)&=\Big(F(z)\E[Y(z)]- \Big(\frac{G^2(z)}{D^2(z)}S(z)\widehat{Y}(z)\nonumber\\
 &+F(z) \int_0^z \frac{G^2(s,a)}{D^2(s,a)}\sum_{n=0}^{\infty}J_n F((s,a);z) S(s,a)\widehat{Y}(s,a)dsda\Big)\Big)dz\nonumber\\
     &+\Big(\frac{G(z)}{D^2(z)}S(z)+F(z) \int_0^z \frac{G(s,a)}{D^2(s,a)}\sum_{n=0}^{\infty}J_n F((s,a);z) S(s,a)dsda\Big)U(dz),
     \end{align*}
     \item
     where the error function 
     $$S(z):= \E[\widetilde{Y}^2(z)]=\E[(Y(z)-\widehat{Y}(z))^2]$$
     satisfies the following time-space Riccati type partial differential equation: 
     \begin{align}
&\frac{\partial ^2}{\partial t \partial x} S(z)
         =2F(z)\Big(\E\big[Y^2(z)\big]-\big(\mathbb{E}[Y(z)]^2\Big) +C^2(z) \nonumber\\
    &+\int_0^t \int_0^x F(s,x)F(t,a) \Big\{\E\big[Y^2(s,a)\big]-2 \E[Y(s,x)]\E[Y(t,a)]\Big\} ds da\nonumber\\  
    &+ \frac{G^2(z)}{D^2(z)} S^2(z)+\Lambda(z); \quad S(t,0)=S(0,x)=\E[(Y(0)-\E[Y(0)])^2],\label{5.11}
    \end{align}
    where
    \small
    \begin{align}
  \Lambda(z)= \frac{\partial ^2}{\partial t \partial x}\Big(\int_0^z \frac{G^2(s,a)}{D^2(s,a)} \Big(2\sum_{n=1}^{\infty}J_n F((s,a);z)+\Big(\sum_{n=1}^{\infty}J_n F((s,a);z)\Big)^2\Big) S^2(s,a) dsda\Big) \label{Lambda}      
    \end{align}
    \end{myenumerate}
\end{theorem}
\dproof
By the Pythagoras theorem we have
\begin{align}
    S(z)=\E[(Y(z)-\widehat{Y}(z))^2] = \E[Y^2(z)]-\E[\widehat{Y}^2(z)]\label{5.12}
\end{align}
From \eqref{2.14a} we have
\begin{align}
\frac{\partial ^2}{\partial t \partial x} \E\Big[Y^2(z)\Big]
=2F(z)\E\Big[Y^2(z)\Big] +C^2(z)+\int_0^t \int_0^x F(s,x)F(t,a) \E\Big[Y^2(s,a)\Big] ds da.\label{5.4}
\end{align}
From \eqref{5.3} we have 
\begin{align}
&\frac{\partial ^2}{\partial t \partial x} \E\Big[\widehat{Y}^2(z)\Big]\nonumber\\
&=\frac{\partial ^2}{\partial t \partial x}\Big(\E\Big[\E[Y(z)] + \int_0^z \Big(\frac{G(s,a)}{D(s,a)} \sum_{n=0}^{\infty}J_n F((s,a);z)S(s,a)\label{5.5}\Big)  M(ds,da)\Big]^2\Big)\nonumber\\
&=\frac{\partial ^2}{\partial t \partial x}\Big(\E[Y(z)]^2 + \int_0^z \frac{G^2(s,a)}{D^2(s,a)} \Big(1+\sum_{n=1}^{\infty}J_n F((s,a);z)\Big)^2 S^2(s,a) dsda\Big)\nonumber\\
&=\frac{\partial ^2}{\partial t \partial x}\Big(\E[Y(z)]^2 \Big) +  \frac{G^2(z)}{D^2(z)}S^2(z)+\Lambda(z),
\end{align}
where 
\begin{align}
  \Lambda(z)= \frac{\partial ^2}{\partial t \partial x}\Big(\int_0^z \frac{G^2(s,a)}{D^2(s,a)} \Big(2\sum_{n=1}^{\infty}J_n F((s,a);z)+\Big(\sum_{n=1}^{\infty}J_n F((s,a);z)\Big)^2\Big) S^2(s,a) dsda\Big).
\end{align}
In general we have
\begin{align}
    \frac{\partial ^2}{\partial t \partial x}(H(z)K(z))&=H(z)\frac{\partial ^2}{\partial t \partial x}K(z)+K(z)\frac{\partial ^2}{\partial t \partial x}H(z)\nonumber\\
    &+\frac{\partial }{\partial x}H(z) \frac{\partial}{\partial t}K(z)+\frac{\partial }{\partial x}K(z) \frac{\partial}{\partial t}H(z).\label{5.16}
\end{align}
Applying this to $H(z)=K(z)=\E[Y(z)]$, we get
\begin{align}
    \frac{\partial ^2}{\partial t \partial x}\Big(\E[Y(z)]^2\Big)&=2F(z) \Big(\E[Y(z)]\Big)^2 \nonumber\\
    &+2\Big(\int_0^xF(t,a) \E[Y(t,a)]da\Big) \Big(\int_0^tF(s,x)\E[Y(s,x)]ds\Big),\label{5.15}
\end{align}
Substituting this in \eqref{5.5} we get
\begin{align}
&\frac{\partial ^2}{\partial t \partial x} \E\Big[\widehat{Y}^2(z)\Big]=\frac{\partial ^2}{\partial t \partial x}\Big(\E[Y(z)]^2 \Big) +  \frac{G^2(z)}{D^2(z)} \Big(\sum_{n=0}^{\infty}J_n F((s,a);z)\big)^n\Big)^2 S^2(z)\nonumber\\
&=2F(z) \Big(\mathbb{E}[Y(z)]\Big)^2 
    +2\Big(\int_0^x F(t,a) \E[Y(t,a)]da\Big) \Big(\int_0^t F(s,x)\E[Y(s,x)]ds\Big)\nonumber\\
    &+ \frac{G^2(z)}{D^2(z)} S^2(z)+\Lambda(z). \label{5.17}
    \end{align}
    Substituting \eqref{5.15} and \eqref{5.17} in \eqref{5.11} we obtain
    \begin{align*}
         \frac{\partial ^2}{\partial t \partial x} S(z)&=\frac{\partial ^2}{\partial t \partial x}\Big( \E[Y^2(z)]\Big)-\frac{\partial ^2}{\partial t \partial x}\Big(\E[\widehat{Y}^2(z)]\Big)\nonumber\\
         &=2F(z)\E\Big[Y^2(z)\Big] +C^2(z)+\int_0^t \int_0^x F(s,x)F(t,a) \E\Big[Y^2(s,a)\Big] ds da\nonumber\\
         &-\Big\{2F(z) \Big(\big(\mathbb{E}[Y(z)]\Big)^2 
    +2\Big(\int_0^x F(t,a) \E[Y(t,a)]da\Big) \Big(\int_0^t F(s,x)\E[Y(s,x)]ds\Big)\nonumber\\
    &+ \frac{G^2(z)}{D^2(z)} S^2(z)+\Lambda(z)\Big\}\nonumber\\
    &=2F(z)\Big(\E\big[Y^2(z)\big]-\big(\mathbb{E}[Y(z)]^2\Big) +C^2(z) \nonumber\\
    &+\int_0^t \int_0^x F(s,x)F(t,a) \Big\{\E\big[Y^2(s,a)\big]-2 \E[Y(s,x)]\E[Y(t,a)]\Big\} ds da\nonumber\\  
    &+ \frac{G^2(z)}{D^2(z)} S^2(z)+\Lambda(z).
    \end{align*}
    \fproof
\begin{remark}
The functions    $\E[Y(z)]$ and $\E[Y^2(z)]$ can be found by \eqref{2.12a} and \eqref{2.14a} respectively.Then $S(z)$ is found by solving the Riccati equation \eqref{5.11}.
Note that these 3 quantities do not depend on the observations and can therefore be computed before the observations start.
\end{remark} 

\section{Examples: Filtering problem with noisy observations}
This section is devoted to discussing some applications of the filtering problem, which involve noisy observations of a given random variable and noisy observations of a Brownian sheet.
\subsection{(Noisy observations of a constant)}
Suppose the signal process is simply a (constant) Gaussian random variable $\theta$, i.e.
\begin{align*}
    dY(z)=0;\quad Y(0,x)=Y(t,0)=\theta,
\end{align*}
and we only have noisy observations of this random variable. Then the observation process is simply
\begin{align}
    dU(z)=\theta dz + B_2(dz); \quad U(0,x)=U(t,0)=\E[\theta].\label{6.2a}
\end{align}
Using Theorem \ref{5.1} with $F=C=0$ and $G=D=1$ we deduce the following:
\begin{theorem}
The best estimate $\widehat{Y}(z)$ of the constant $\theta$ given the noisy observations \eqref{6.2a} solves the SDE
    \begin{align*}
        d\widehat{Y}(z)= -S(z)\widehat{Y}(z) dz + S(z) U(dz); \quad \widehat{Y}(0,x)=\widehat{Y}(t,0)=\E[\theta],
    \end{align*}
    where $S(z)$ satisfies the time-space Riccati equation
 \begin{align*}
     dS(z)=S^2(z) dz; \quad S(0,x)=S(t,0)=\E[(\theta - \E[\theta])^2].
 \end{align*}   
\end{theorem}

\subsection{(Noisy observations of the Brownian sheet)}
Suppose the signal process is the Brownian sheet, i.e. 
\begin{align*}
    dY(z)=B_1(dz); Y(0,x)=Y(t,0)=0,
\end{align*}
and we make noisy observations of this process. Then the observation process is
\begin{align}
    dU(z)= Y(z)dz + B_2(dz); U(0,x)=U(0,t)=0.\label{6.2}
\end{align}
Then by using  Theorem \ref{5.1} with $F=0$ and $C=D=G=1$ we get the following result:
\begin{theorem} The best estimate $\widehat{Y}(z)$ of $Y(z)=B_1(z)$ given the noisy observations \eqref{6.2} satisfies the following SDE
\begin{align*}
 d\widehat{Y}(z)&= -S(z)\widehat{Y}(z)dz +  S(z) U(dz);\quad \widehat{Y}(0,x)=\widehat{Y}(t,0)=0,      
\end{align*}
where $S(z)$ satisfies the time-space Riccati equation
     \begin{align*}
dS(z)       = (1 + S^2(z))dz; \quad S(0,x)=S(t,0)=0.\label{6.4}
    \end{align*}
\end{theorem} 
\begin{remark}
    We remark that the time-space Riccati equation always has a solution (see Remark 7.5 \cite{AOPT1}).
\end{remark}

\textbf{Acknowledgements}
Nacira Agram and and Olena Tymoshenko gratefully acknowledge the financial support provided by the Swedish Research Council grants (2020-04697), the Slovenian Research and
Innovation Agency, research core funding No.P1-0448, and the MSCA4Ukraine grant (AvH ID:1233636), respectively.

\end{document}